\documentclass[a4paper, leqno, 11pt]{article}
\usepackage{latexsym}
\usepackage{amsmath}
\usepackage{amssymb}
\usepackage{enumerate}
\usepackage{theorem}
\usepackage{array}
\newtheorem{itheorem}{Theorem}
\newtheorem{theorem}{Theorem}[section]
\newtheorem{lemma}[theorem]{Lemma}

\newtheorem{proposition}[theorem]{Proposition}
\theorembodyfont{\rmfamily}

\newtheorem{setup}[theorem]{Setup}
\newtheorem{remark}[theorem]{Remark}

\newtheorem{conjecture}{Modified Sziklai's Conjecture}

\newcommand{\proof}{\noindent \mbox{\em Proof.\hspace*{2mm}}}
\newcommand{\qed}{\hfill \mbox{$  \Box $}}
\newcommand{\gyokan}{\vskip 4pt}
\title{Sziklai's conjecture on the number of points of
a plane curve over a finite field II
\footnote{
The final version appeared in Contemporary Mathematics 518 (2010), 225--234.
This updated edition contains an addendum,
which fixes an error in Remark~\ref{incorrectremark} of the printed version.}
}
\author{Masaaki Homma
\thanks{Partially supported by Grant-in-Aid
for Scientific Research (21540051), JSPS.}
\\
 Department of Mathematics,
Kanagawa University\\
Yokohama 221-8686, Japan\\
homma@n.kanagawa-u.ac.jp
\and Seon Jeong Kim
\thanks{Partially supported by 
the Korea Research Foundation Grant funded by the
Korean Government(MOEHRD) (KRF-2006-312-C00016).}\\
Department of Mathematics and RINS\\
Gyeongsang National University\\
Jinju 660-701, Korea \\
skim@gnu.kr}
\date{}

\begin{document}

\maketitle
\begin{abstract}
We settle the conjecture posed by Sziklai
on the number of points of a plane curve over a finite field
under the assumption that the curve is nonsingular.

\end{abstract}

\section{Introduction}
In the paper \cite{szi},
Sziklai posed a conjecture on the number of points of a plane curve over a finite field.
Let $C$ be a plane curve of degree $d$ over ${\Bbb F}_q$
without an ${\Bbb F}_q$-linear component.
Then he conjectured that the number of ${\Bbb F}_q$-points $N_q(C)$ of $C$
would be at most $(d-1)q + 1$.
But he had overlooked the known example of a curve of degree $4$ over ${\Bbb F}_4$ with $14$ points
(\cite{ser}, \cite{gee-vlu}).
So we must modify this conjecture.
\begin{conjecture}
Unless $C$ is a curve defined over ${\Bbb F}_4$
which is projectively equivalent to
{\small
\begin{equation}\label{counterexample}
X^4 + Y^4 +Z^4 +X^2Y^2 +Y^2Z^2 + Z^2X^2
    + X^2YZ +XY^2Z + XYZ^2 =0
\end{equation}
}
over ${\Bbb F}_4$,
we might have 
\begin{equation}\label{guess}
N_q(C) \leq (d-1)q +1.
\end{equation}
\end{conjecture}
Here we make two parenthetical remarks on this conjecture.
Since $C$ is defined by a homogeneous equation
$F(X, Y, Z) = 0$,
we understand the set of ${\Bbb F}_q$-points $C({\Bbb F}_q)$
of $C$ to be the set of ${\Bbb F}_q$-points
$(\alpha , \beta , \gamma ) \in {\Bbb P}^2$
such that $F(\alpha , \beta , \gamma) = 0$, that is to say,
it is no matter whether each of those points is nonsingular or not.
The second remark is that the conjecture makes sense only if
$2 \leq d \leq q+1$ because the conjectural bound exceeds the obvious bound
$N_q(C) \leq {}^{\#}{\Bbb P}^2({\Bbb F}_q) = q^2 + q +1$
if $d \geq q+2$.

In the previous paper \cite{hom-kim2},
we proved the inequality
\begin{equation}\label{previousbound}
N_q(C) \leq d(q-1) +2 = (d-1)q + (q+2-d),
\end{equation}
which guarantees the inequality (\ref{guess}) for $d=q+1$,
and presented an example of a curve of degree $q+1$ having $q^2+1$
${\Bbb F}_q$-points.
Moreover, we observed that if a curve of degree $4$ over ${\Bbb F}_4$
has more than $13$ rational points,
then this curve is projectively equivalent to
the curve (\ref{counterexample}) over ${\Bbb F}_4$.

\gyokan

The main purpose of this paper is to show the following.
\begin{itheorem}
For $d=q$, the modified Sziklai's conjecture holds true,
and for each $q$ there exists a nonsingular curve of degree $q$
over ${\Bbb F}_q$ with $(q -1)q +1$ rational points.
\end{itheorem}
Note that the truth of the inequality (\ref{guess}) for $d=q=3$
is classical \cite{seg}, and it is well known for $d=q=2$.
Additionally, we show that (\ref{guess}) holds if the curve $C$ is
nonsingular of degree $d \leq q-1$.
Therefore, together with our previous results,
the following theorem is established.
\begin{itheorem}
The modified Sziklai's conjecture is true for nonsingular curves.
Moreover there is an example of a nonsingular curve for which
equality holds in \mbox{\rm (\ref{guess})}
if $ d= q+2, q+1, q, q-1, \sqrt{q} +1  
\mbox{{\rm  \ (}when $q$ is square{\rm )}, or } 2$.
\end{itheorem}
\section{Simplification of the problem}\label{Sec_sim}
To settle the modified Sziklai's conjecture affirmatively,
{\em we may suppose the curve $C$ to be absolutely irreducible without
an ${\Bbb F}_q$-rational singular point}.
Actually the following three facts hold.
Throughout this section, we assume that the degree of $C$ is at most $q+1$.
\begin{proposition}\label{irredfq}
If $C$ is reducible over ${\Bbb F}_q$,
then $N_q(C) < (d-1)q$.
\end{proposition}
\begin{proposition}\label{noirredcompnotfq}
If $C$ has an irreducible component which is not defined over ${\Bbb F}_q$,
then $N_q(C) \leq (d-1)q$.
\end{proposition}
\begin{proposition}\label{nofqsing}
If $C$ has a singular point which is an ${\Bbb F}_q$-point,
then $N_q(C) \leq (d-1)q$.
\end{proposition}
\noindent
{\em Proof of Proposition}~\ref{irredfq}.\hspace*{2mm}
Let $C = C_1 \cup C_2$,
where each curve $C_i$ is of degree $d_i$, and defined over ${\Bbb F}_q$
without an ${\Bbb F}_q$-linear component.
By a theorem of Segre~\cite[Teorema II on page 30]{seg}
\[
N_q(C_i) \leq (d_i-1)q + \left\lfloor \frac{d_i}{2} \right\rfloor
 \hspace{10mm}
 (i= 1, 2),
\]
where $\lfloor \frac{d_i}{2} \rfloor$ denotes the integer part of
$\frac{d_i}{2}$.
Hence $N_q(C)\leq N_q(C_1) + N_q(C_2) \leq (d-1)q$
because $d_1 + d_2 = d \leq q+1$.
\qed
\gyokan
\noindent
{\em Proof of Proposition}~\ref{noirredcompnotfq}.\hspace*{2mm}
Let $C_1$ be an irreducible component of $C$
which is not defined over ${\Bbb F}_q$,
and ${\Bbb F}_{q^t}$ the minimum extension of ${\Bbb F}_q$
over which $C_1$ is defined.
Since the $t$ conjugates $C_1, \ldots , C_t$ of $C_1$ over ${\Bbb F}_q$
are components of $C$,
$C = C' \cup C_1 \cup \ldots \cup C_t$,
where $C'$ is a curve defined over ${\Bbb F}_q$ or $C' = \emptyset$.
Let $e = \deg C_1$, so $\deg C_1 = \ldots = \deg C_t = e$.
Since $(C_1\cup \ldots \cup C_t)({\Bbb F}_q)
\subset C_1\cap \ldots \cap C_t$,
$N_q(C_1\cup \ldots \cup C_t) \leq e^2$ by B\'{e}zout's theorem.

When $C' \neq \emptyset$,
it is a case of Proposition~\ref{irredfq}.
So we may suppose $C' = \emptyset$.
Then $d = te$ and $N_q(C) \leq e^2$.
Since
\[
t \left( (d-1)q - e^2 \right)
\geq 2(d-1)q -de = d(q-e) + (d-2)q \geq 0,
\]
we have $e^2 \leq (d-1)q$.
\qed

\gyokan
\noindent
{\em Proof of Proposition}~\ref{nofqsing}.\hspace*{2mm}
Let $P_0$ be a singular and ${\Bbb F}_q$-rational point of $C$.
Then for each ${\Bbb F}_q$-line $l$ passing through $P_0$,
${}^{\#}(l \setminus \{P_0\}) \cap C({\Bbb F}_q) \leq d-2$.
So $N_q(C) \leq (d-2)(q+1) +1 \leq (d-1)q$
because $d \leq q+1$.
\qed
\section{The proof for the case $d=q >4$}
Throughout this section,
we fix a plane curve $C$ over ${\Bbb F}_q$ of degree $q$
without an ${\Bbb F}_q$-linear component.
Suppose that $C({\Bbb F}_q) \neq \emptyset$.
\begin{proposition}\label{blackpoint}
Fix an ${\Bbb F}_q$-point $P_0 \in C$,
and an ${\Bbb F}_q$-line $l_{\infty} \subset {\Bbb P}^2$
with $l_{\infty} \not\ni P_0$.
Suppose there are ${\Bbb F}_q$-lines $l_1, \ldots , l_t$
with $q \geq t \geq 3$ passing through $P_0$ such that
the $q$ ${\Bbb F}_q$-points of $l_i \setminus l_{\infty}$
are contained in $C$.
For an ${\Bbb F}_q$-line $l \ni P_0$
other than these $t$ lines,
if ${}^{\#} \left( (l\setminus l_{\infty}) \cap C({\Bbb F}_q) \right)
\geq q-t+2$,
then all the $q$ ${\Bbb F}_q$-points of $l\setminus l_{\infty}$
are contained in $C$.
\end{proposition}
\proof
Choose coordinates $X, Y, Z$ of ${\Bbb P}^2$ as
$l_1$ is defined by $X=0$, $l_2$ by $Y=0$, and $l_{\infty}$ by $Z=0$.
So $P_0 =(0,0,1)$.
Let
\[
f(x,y) = \sum_{i, j {\rm \ with}
               \atop
               i+j \leq q}
               a_{ij} x^i y^j = 0
\]
be an affine equation over ${\Bbb F}_q$ defining $C$ on the affine plane
${\Bbb P}^2 \setminus l_{\infty}$ with affine coordinates
$x = \frac{X}{Z}$, $y= \frac{Y}{Z}$.
Since $l_1({\Bbb F}_q) \subset C$, $f(0, \beta) =0$ for any
$\beta \in {\Bbb F}_q$.
Hence $f(0, y) = a_{0q}(y^q-y)$
because the degree of $f(0, y)$ is at most $q$.
Similarly, $f(x, 0) = a_{q0}(x^q-x).$
Hence
\[
f(x,y) = a_{q0}(x^q-x) + a_{0q}(y^q-y) + 
    xy\left( g_{q-2}(x,y) + \ldots + g_1(x,y) + g_0 \right),
\]
where
$g_{\nu}(x,y) = \sum_{k=0}^{\nu} a_{\nu - k+1, k+1} x^{\nu-k}y^k$.
Let $y = u_{\mu}x $ ($u_{\mu} \in {\Bbb F}_q^{\times}$) be an affine
equation of the line $l_{\mu}$ for $\mu = 3, \ldots, t$.
Here ${\Bbb F}_q^{\times}$ denotes the multiplicative group
${\Bbb F}_q \setminus \{ 0 \}.$
Since $f(\alpha , u_{\mu}\alpha) =0$ for any $\alpha \in {\Bbb F}_q$
by the assumption on $l_{\mu}$,
we have
\[
\left(
  \begin{array}{lllll}
    & \vdots &   & & \\
   \alpha^{q-2}&\alpha^{q-3}& \cdots &\alpha & 1 \\
    & \vdots &   & &
  \end{array}
\right)_{\alpha \in {\Bbb F}_q^{\times}}
 \left(
  \begin{array}{c}
   g_{q-2}(1,u_{\mu})\\
   g_{q-3}(1,u_{\mu})\\
    \vdots \\
   g_1(1,u_{\mu})\\
    g_0
  \end{array}
 \right)
 =
 \left(
  \begin{array}{c}
   0\\
   \vdots \\
   0
  \end{array}
 \right).
\]
Since $\det 
\left( \alpha^{k} \right)_{
    (\alpha, k) \in {\Bbb F}_q^{\times} \times \{q-2, \ldots , 0\}}
 \neq 0$,
\[
g_{q-2}(1,u_{\mu}) = \ldots = g_1(1,u_{\mu}) = g_0 = 0.
\]
Hence each equation $g_{\nu}(1,y)=0$ has at least $t-2$ zeros,
which implies
$g_{t-3}(1,y) = \ldots = g_0 = 0$ as polynomials.
So
\[
f(x,y) = a_{q0}(x^q-x) + a_{0q}(y^q-y) + 
    xy \sum_{\nu = t-2}^{q-2} g_{\nu}(x,y).
\]
Let $y =vx$ ($v \in {\Bbb F}_q^{\times}$) be an equation of $l$.
By the assumption on $l$, there are at least $q-t+1$ elements
$\alpha_1, \ldots, \alpha_{q-t+1} \in {\Bbb F}_q^{\times}$
so that $f(\alpha_i, v\alpha_i) = 0$
for each $i = 1, \ldots, q-t+1$.
Hence
\[
 \left(
  \begin{array}{cccc}
   & \vdots &  &  \\
  \alpha_i^{q-2}&\alpha_i^{q-3}& \cdots & \alpha_i^{t-2}\\
   & \vdots &  &  
  \end{array}
 \right)_{i = 1, \ldots , q-t+1}
 \left(
   \begin{array}{c}
     g_{q-2}(1,v) \\
     \vdots \\
     g_{t-2}(1,v)
   \end{array}
 \right)
=
 \left(
  \begin{array}{c}
   0\\
   \vdots \\
   0
  \end{array}
 \right).
\]
Since
$\det \left( \alpha_i^k \right) \neq 0$,
we have
$g_{q-2}(1,v) = \ldots = g_{t-2}(1,v) =0$, and get
$f(1,v) =0$.
\qed

\begin{proposition}\label{whitepoint}
Fix an ${\Bbb F}_q$-point $Q_0 \in {\Bbb P}^2({\Bbb F}_q) \setminus C.$
Suppose there are ${\Bbb F}_q$-lines $l_1, \ldots , l_t$ with
$q-1 \geq t \geq 2$ passing through $Q_0$ such that
$l_i({\Bbb F}_q) \setminus \{ Q_0 \} \subset C$.
If an ${\Bbb F}_q$-line $l \ni Q_0$ other than these $t$ lines
has at least $q-t+1$ ${\Bbb F}_q$-points of $C$,
then $l({\Bbb F}_q) \setminus \{ Q_0 \} \subset C$
\end{proposition}
\proof
First choose $q-t+1$ points in
$\left( l({\Bbb F}_q) \setminus \{ Q_0 \} \right) \cap C$,
and then choose an ${\Bbb F}_q$-point $P'$ of $l \setminus \{ Q_0 \}$
other than these $q-t+1$ points.
Fix an ${\Bbb F}_q$-line, say $l_{\infty}$,
such that $l_{\infty} \ni P'$ but $l_{\infty} \not\ni Q_0$.
Choose coordinates $X, Y, Z$ of ${\Bbb P}^2$
so that $l_1$ is defined by $X=0$, $l_2$ by $Y=0$,
and $l_{\infty}$ by $Z=0$.
Then $Q_0 =(0,0,1)$.
Let
\[
F(X,Y,Z) = \sum_{i, j {\rm \ with}
               \atop
               i+j \leq q}
               a_{ij}X^iY^jZ^{q-i-j}
\]
be a homogeneous equation over ${\Bbb F}_q$ defining $C$.
Since $l_1({\Bbb F}_q) \setminus \{ Q_0\} \subset C,$
$0 = F(0,1, \beta) = \sum_{j=0}^q a_{0j}\beta^{q-j}$
for any $\beta \in {\Bbb F}_q$.
So $F(0,1,Z) = a_{00}(Z^q-Z)$, and hence
$F(0,Y,Z) = a_{00}(Z^q-Y^{q-1}Z)$.
Similarly $F(X,0,Z) = a_{00}(Z^q-X^{q-1}Z)$.
Therefore
\begin{eqnarray*}
\lefteqn{F(X,Y,Z) = a_{00}(Z^q-X^{q-1}Z-Y^{q-1}Z) +}\\
&&  XY(g_{q-2}(X,Y) + g_{q-3}(X,Y)Z
    + \ldots + g_0 Z^{q-2}),
\end{eqnarray*}
where
$g_{\nu}(X,Y) = \sum_{k=0}^{\nu} a_{\nu-k+1, k+1}X^{\nu-k}Y^{k}.$

In general, any line $L$ over ${\Bbb F}_q$ which contains $Q_0$
but is not $l_1$ nor $l_2$ is defined by an equation of the form
$Y=uX$ for some $u \in {\Bbb F}_q^{\times}$.
So $L({\Bbb F}_q)\setminus \{Q_0\} = \{ (1,u,\beta) | \beta \in {\Bbb F}_q \}$.
Note that
\begin{eqnarray*}
\lefteqn{F(1,u,\beta) = ug_{q-2}(1,u) +(u g_{q-3}(1,u)-a_{00})\beta +}\\
&&  ug_{q-4}(1,u)\beta^2 + \ldots +ug_{1}(1,u)\beta^{q-3} + ug_0 \beta^{q-2},
\end{eqnarray*}
because $\beta^q - \beta - u^{q-1}\beta = - \beta.$
Let $Y=u_{\mu}X$ be an equation of $l_{2 + \mu}$ ($\mu = 1, \ldots, t-2$).
Note that these $u_{\mu}$'s are not $0$.
Then
\[
 \left(
  \begin{array}{ccccc}
    & & \vdots &  &  \\
  1&\beta&\beta^2 & \cdots & \beta^{q-2}\\
    & & \vdots &  &  
  \end{array}
 \right)_{\beta \in {\Bbb F}_q^{\times}}
 \left(
   \begin{array}{c}
     u_{\mu}g_{q-2}(1, u_{\mu})\\
      u_{\mu}g_{q-3}(1, u_{\mu})-a_{00}\\
       u_{\mu}g_{q-4}(1, u_{\mu})\\
     \vdots \\
     u_{\mu}g_{0}
   \end{array}
 \right)
=
 \left(
  \begin{array}{c}
   0\\
   \vdots \\
   0
  \end{array}
 \right).
\]
Since
$\det \left(
       \beta^{k}
      \right)_{(\beta, k) \in 
         {\Bbb F}_q^{\times} \times \{0,1, \ldots, q-2\}}
          \neq 0,$
we have, in particular,
$
u_{\mu}g_{q-4}(1, u_{\mu}) = \ldots = u_{\mu}g_{0} =0.
$
Hence if $\nu < t-2$, then $g_{\nu}(1,y) =0$
as a polynomial in $y$, 
because $g_{\nu}(1,y) =0$ has $t-2$ roots
$\{ u_3, u_4, \ldots , u_t \}$
but its degree is less than $t-2$.
Therefore
\begin{eqnarray*}
\lefteqn{F(1,y,z) = a_{00}(z^q-y^{q-1}z) +}\\
&&  yg_{q-2}(1,y) + (yg_{q-3}(1,y) -a_{00})z
    +yg_{q-4}(1,y)z^2+ \ldots + yg_{t-2}(1,y) z^{q-t}.
\end{eqnarray*}
Let $Y=vX$ be an equation of $l$,
and
$\{ (1, v, \beta_i) | 1\leq i \leq q-t+1 \}$
a set of chosen points of
$(l({\Bbb F}_q) \setminus \{Q_0\}) \cap C.$
Then
\[
 \left(
  \begin{array}{ccccc}
    & & \vdots &  &  \\
  1&\beta_i&\beta_i^2 & \cdots & \beta_i^{q-t}\\
    & & \vdots &  &  
  \end{array}
 \right)_{i = 1, \ldots, q-t+1}
 \left(
   \begin{array}{c}
     vg_{q-2}(1, v)\\
      vg_{q-3}(1, v)-a_{00}\\
       vg_{q-4}(1, v)\\
     \vdots \\
     vg_{t-2}(1,v)
   \end{array}
 \right)
=
 \left(
  \begin{array}{c}
   0\\
   \vdots \\
   0
  \end{array}
 \right).
\]
Hence
$vg_{q-2}(1, v) = vg_{q-3}(1, v)-a_{00} = \ldots = vg_{t-2}(1,v) =0,$
and then $F(1,v,\beta)=0$ for any $\beta \in {\Bbb F}_q$,
which means that $l({\Bbb F}_q) \setminus \{ Q_0 \} \subset C$.
\qed

\gyokan

Now we prove the following theorem by a reduction to absurdity.
\begin{theorem}\label{theorem_q}
Let $C$ be a plane curve over ${\Bbb F}_q$ of degree $q$
without an ${\Bbb F}_q$-linear component.
If $q>4$, then $N_q(C) \leq (q-1)q +1$.
\end{theorem}
By the previous result (\ref{previousbound}),
$N_q(C) \leq (q-1)q +2.$
We prove the absurdity of the equality
$N_q(C) = (q-1)q +2.$
Moreover, by the arguments in Section~\ref{Sec_sim},
we may assume that $C$ is irreducible and each ${\Bbb F}_q$-rational
point of $C$ is nonsingular.
\begin{setup}\label{setup_q}
Until the end of this section,
we suppose that $C$ is an irreducible plane curve of degree $q$
over ${\Bbb F}_q$ with $N_q(C)= (q-1)q +2$ and no point of $C({\Bbb F}_q)$
is singular.
\end{setup}
Some symbols should be introduced here.
$\check{\Bbb P}^2$ is the projective plane of lines
in the original plane ${\Bbb P}^2$.
So $\check{\Bbb P}^2({\Bbb F}_q)$ means
the set of ${\Bbb F}_q$-lines of ${\Bbb P}^2$.
Let
\[
{\cal A}_i = \{ l \in \check{\Bbb P}^2({\Bbb F}_q) \mid
     {}^{\#}\left( l \cap C({\Bbb F}_q) \right) = i \}
\]
and $a_i = {}^{\#}{\cal A}_i.$
\begin{lemma}\label{numberofpointonl}
Under Setup~{\rm \ref{setup_q}}, we have
\begin{enumerate}[{\rm (1)}]
\item $\displaystyle 
        \sum_{i=0}^q a_i = q^2 + q +1;$
\item $\displaystyle
        \sum_{i=0}^q  ia_i = (q+1)(q^2 -q +2);$
\item $\displaystyle
        \sum_{i=2}^q \left( \begin{array}{c}
                      i\\
                      2
                    \end{array}\right)  a_i
                 = \left( \begin{array}{c}
                      q^2-q +2\\
                      2
                    \end{array}\right);$
\item $\displaystyle
         \sum_{i=1}^{\lfloor \frac{q}{2} \rfloor}ia_i
             + \sum_{j=1}^{\lfloor \frac{q-1}{2} \rfloor}ja_{q-j}
                      \geq q^2 -q +2.$
\end{enumerate}
\end{lemma}
\proof
(1) is obvious.
For (2), consider the point-line correspondence
\[
 {\cal P} = \{ (P,l) \in C({\Bbb F}_q)\times \check{\Bbb P}^2({\Bbb F}_q)
               \mid P \in l \}
\]
with two projections $\pi_1: {\cal P} \to C({\Bbb F}_q)$ and
$\pi_2: {\cal P} \to \check{\Bbb P}^2({\Bbb F}_q)$.
Counting the number ${}^{\#}{\cal P}$ by using $\pi_1$,
we have $(q+1)(q^2 -q +2)$, and by $\pi_2$,
$\sum_{i=0}^q \,  {}^{\#}  \pi_2^{-1}({\cal A}_i ) = \sum_{i=0}^q  ia_i.$
  
For (3), consider the correspondence
\[
 {\cal P}'' = \{ (\{P, Q\} , l ) \in 
      \left( S^2C({\Bbb F}_q) \setminus \Delta \right) \times
               \check{\Bbb P}^2({\Bbb F}_q)
            \mid P, Q \in l \},
\]
where $S^2C({\Bbb F}_q)$ denotes the symmetric product of
two copies of $C({\Bbb F}_q)$ and $\Delta$ the diagonal subset of
$S^2C({\Bbb F}_q)$.
Counting ${}^{\#}{\cal P}''$ by using the first projection
$\pi_1'': {\cal P}'' \to S^2C({\Bbb F}_q) \setminus \Delta$
and the second projection
$\displaystyle \pi_2'': {\cal P}'' \to \check{\Bbb P}^2({\Bbb F}_q) 
   = \coprod_{i=0}^q {\cal A}_i$,
we have the desired formula.
  
For (4), consider the correspondence
\[
 {\cal P}' = \{ (P,l) \in C({\Bbb F}_q)\times \check{\Bbb P}^2({\Bbb F}_q)
               \mid i(l.C;P) \geq 2 \},
\]
where $i(l.C;P)$ denotes the intersection multiplicity of $l$ and $C$
at $P$.
Note that for each point $P \in C({\Bbb F}_q)$,
there is a unique ${\Bbb F}_q$-line $l$ such that $i(l.C;P) \geq 2$
because $C$ is nonsingular at $P$.
So ${}^{\#}{\cal P}' = q^2 -q +2.$
Let $\pi_2' : {\cal P}' \to \check{\Bbb P}^2({\Bbb F}_q)$ be
the second projection.
For a line $l \in {\cal A}_i$,
let ${}^{\#}\left({\pi_2'}^{-1}(l)\right) = s_l.$
For each point $P$ of these $s_l$ points on $l$,
$i(l.C;P) \geq 2$ by definition.
Hence we have $2s_l + i -s_l \leq (l.C) = q$
by B\'{e}zout's theorem.
So $s_l \leq \min \{i, q-i\}.$
Hence ${}^{\#}{\cal P}' \leq \sum_{i=1}^{q} \min \{i, q-i\} a_i.$
\qed

\gyokan

\noindent
{\em Proof of Theorem}~\ref{theorem_q}.\hspace*{2mm}
\underline{Step~I}. We prove that $a_0 = a_1 = 0$.

By (1), (2) and (4) of Lemma~\ref{numberofpointonl},
\begin{eqnarray*}
 \lefteqn{qa_0 + 
   \sum_{i=1}^{\lfloor \frac{q}{2} \rfloor} (q-2i)a_i }\\
  &=& q(\sum_{i=0}^{q} a_i) - \sum_{i=0}^{q} ia_i - 
        (\sum_{i=1}^{\lfloor \frac{q}{2} \rfloor} ia_i
           + \sum_{j=1}^{\lfloor \frac{q-1}{2} \rfloor} ja_{q-j} ) \\
  & \leq & q(q^2+ q+1) - (q+1)(q^2-q+2) - (q^2-q +2) \\
  &= & q-4.
\end{eqnarray*}
Hence $a_0 = a_1 = 0.$

\noindent
\underline{Step~II}.
We prove that $a_2=0$.

Suppose $a_2 > 0$. Choose a line $l_0 \in {\cal A}_2.$
Two of the $q+1$ ${\Bbb F}_q$-points of $l_0$ are on $C$,
say $P_0$ and $P_1$, and the other $q-1$ are not on $C$,
say $P_2', \ldots , P_q'.$
Let $l_0, l_1, \ldots , l_q$ be the set of ${\Bbb F}_q$-lines
passing through $P_0$.
For each line $l_i$ with $ 1 \leq i \leq q$,
there is an ${\Bbb F}_q$-point $Q_i$ of $l_i$ not lying on $C$
because ${}^{\#} (l_i \cap C ) \leq q$.
Note that ${\Bbb P}^2({\Bbb F}_q) \setminus C({\Bbb F}_q) =
\{ P_2, \ldots , P_q', Q_1, \ldots , Q_q \}$
because $N_q(C) = (q-1)q +2$.
So $Q_i$ is the unique ${\Bbb F}_q$-point of $l_i$
which does not lie on $C$.
If one considers the all lines passing through $P_1$,
say $l_0, l_1', \ldots l_q'$,
each line $l_i'$ ($ 1 \leq i \leq q$) has a unique ${\Bbb F}_q$-point
not lying on $C$.
So we may assume $l_i' \ni Q_i$ for $i = 1, \ldots, q$.
Hence the line $Q_iQ_j$ never meets with $P_0$ nor $P_1$.
In particular, $P_0, Q_i, Q_j$ are not collinear
for any $1 \leq i, j \leq q$ with $i \neq j$,
and neither are $P_1, Q_i, Q_j$.
If three of $\{ Q_1, \ldots , Q_q \}$ are collinear,
so are $Q_1, \ldots , Q_q$ by Proposition~\ref{blackpoint},
which is a contradiction by Step~I.
Therefore
${\cal K} = \{ P_0, P_1, Q_1, \ldots , Q_q \}$ forms a $(q+2)$-arc.
Hence $q$ must be a power of $2$ \cite[Theorem~8.5]{hir}.
So $q \geq 8$ because $q>4$ a priori.

Next let us consider the $q$ ${\Bbb F}_q$-lines passing through
$P_2'$ other than $l_0$, say $m_1, \ldots , m_q$.
It is easy to see that a $(q+2)$-arc has no unisecant.
So half of $m_1, \ldots , m_q$ do not meet ${\cal K}$,
and each line of the other half meets ${\cal K}$ at two points.
Since $q \geq 8$, we may assume that $m_i \cap {\cal K} = \emptyset$
($i=1,2,3$) and ${}^{\#}(m_4 \cap {\cal K}) = 2.$
Applying Proposition~\ref{whitepoint}
to $P_2', m_1, m_2, m_3$ and $m_4$ as $l$,
we have $m_4$ is also an external line to ${\cal K}$,
which is a contradiction.

\noindent
\underline{Step~III}. 
Let $k = \min \{ i \mid {\cal A}_i \neq \emptyset\}$.
We prove that $3 \leq k \leq q-3$.

We already saw $k \geq 3$.
Suppose that $k \geq q-2$, namely,
$a_0 = a_1 = \ldots = a_{q-3} =0.$
Hence, by Lemma~\ref{numberofpointonl}, we have
\begin{eqnarray}
&&\label{III1} a_{q-2} + a_{q-1} + a_{q} = q^2 + q +1 \\
&&\label{III2} (q-2)a_{q-2} +(q-1) a_{q-1} + qa_{q}
             = (q+1)(q^2 -q +2) \\
&&\label{III3} \left(\begin{array}{c}
     q-2 \\
     2
   \end{array}\right) a_{q-2}
   +
   \left(\begin{array}{c}
     q-1 \\
     2
   \end{array}\right) a_{q-1}
   +
   \left(\begin{array}{c}
     q \\
     2
   \end{array}\right) a_{q}
      =
   \left(\begin{array}{c}
     q^2 -q + 2 \\
     2
   \end{array}\right) .
\end{eqnarray}
Making $-q(q-2)$ times Eq.~(\ref{III1})
plus $2q-3$ times Eq.~(\ref{III2})
minus $2$ times Eq.~(\ref{III3}),
we know
$a_{q-1} = (q-2)(4-q)$,
which is impossible because $q >4$.

\noindent
\underline{Step~IV}.
Fix a line $l_0 \in {\cal A}_k$,
where $k$ is the number explained in the previous step.
Let $P_0, \ldots , P_{k-1}$ be
the $k$ ${\Bbb F}_q$-points of $l_0$ that lie on $C$,
and $P_k', \ldots , P_q'$
the remaining ${\Bbb F}_q$-points of $l_0$.
Let $S = \left( {\Bbb P}^2({\Bbb F}_q)  \setminus C \right)
      \setminus \{ P_k', \ldots , P_q' \}$.
Since ${}^{\#}\left( {\Bbb P}^2({\Bbb F}_q)  \setminus C \right)
     = 2q-1$,  ${}^{\#} S = q+k-2$.
In this step, we show that there is a point $Q \in S$ so that
${}^{\#} \{ P_iQ \mid 0 \leq i \leq k-1, \, P_iQ \in {\cal A}_q \}
    \geq 3.$

Consider the correspondence
\[
{\cal A}' = \{ (Q, P_i) \mid Q \in S, \, i = 0, \ldots , k-1, \,
           P_iQ \in {\cal A}_q \}
\]
with projections
$\pi_1 : {\cal A}' \to S$
and $\pi_2 : {\cal A}' \to \{ P_0, \ldots , P_{k-1} \}.$
Our claim is that there is a point $Q \in S$
so that ${}^{\#}\pi_1^{-1}(Q) \geq 3$.
For each line $l \ni P_i$ except $l_0$, $l \cap S \neq \emptyset$
because $\deg C = q$,
and these $q$ lines $\ni P_i$ cover $S$.
Hence we may suppose that $S = \{ Q_1, Q_2, \ldots , Q_{q+k-2} \}$
and those $q$ lines are $P_iQ_1, \ldots , P_iQ_q$.
Under this notation,
$P_iQ_j \in {\cal A}_q$ if and only if
$P_iQ_j \cap \{ Q_{q+1}, \ldots , Q_{q+k-2} \} = \emptyset$.
So we have ${}^{\#} \pi_2^{-1} (P_i) \geq q - (k-2)$,
and then ${}^{\#}{\cal A}'\geq k(q-k+2)$.
If ${}^{\#}\pi_1^{-1}(Q) \leq 2$ for any $Q \in S$,
we have
${}^{\#}{\cal A}' \leq 2(q+k-2)$.
So $2(q+k-2) - k(q-k+2)$ must be nonnegative.
But this number is equal to
$(k-2)(k - (q-2))$,
which is a contradiction
because $3 \leq k \leq q-3$ by Step~III.

\noindent
\underline{Step~V}. 
Choose a point $Q \in S$ having the property described in the previous step.
We may suppose that $P_iQ \in {\cal A}_q$ for $ i = 0, 1, \ldots , s-1$
with $3 \leq s \leq k$,
and other $q+1-s$ lines $P_sQ, \ldots , P_{k-1}Q, P_{k}'Q, \ldots , P_{q}'Q$
passing through $Q$ do not belong to ${\cal A}_q$.
Let $m$ be one of these $q+1-s$ lines.
Then ${}^{\#} (m \cap C({\Bbb F}_q)) \leq q-s$, otherwise
$m \in {\cal A}_q$ by Proposition~\ref{whitepoint}.
Hence
\[
{}^{\#}(m \cap (S \setminus \{ Q \}))
  \geq \left\{
        \begin{array}{cl}
           s & \mbox{\rm if  $ m= P_iQ$ ($s \leq i \leq k-1$)}\\
          s-1& \mbox{\rm if  $ m= P_j'Q$ ($k \leq j \leq q$).}
        \end{array}
      \right.
\]
Therefore
\[
{}^{\#}(S \setminus \{ Q \})
    \geq s(k-s) + (s-1)(q-k+1).
\]
On the other hand,
since ${}^{\#}S = q+k-2$,
\[
s(k-s) + (s-1)(q-k+1) - {}^{\#}(S \setminus \{ Q \}) = (s-2)(q-1-s) >0,
\]
which is a contradiction.
This completes the proof.
\qed
\section{Nonsingular plane curves of degree at most $q-1$}
In this section, we consider a nonsingular plane curve $C$
over ${\Bbb F}_q$ of degree $d$ with $1 < d \leq q-1$.
\begin{theorem}\label{the_d<q}
Under the above setting,
we have $N_q(C) \leq (d-1)q +1$.
\end{theorem}
\proof
To show this bound,
we need some results by the Brazilian school of curve theory.
We explain those briefly only for plane nonsingular curves.
A nonsingular plane curve $C$ defined over ${\Bbb F}_q$ is said to be
$q$-Frobenius nonclassical if $F_q(P) \in T_P(C)$ for a general 
$\overline{\Bbb F}_q$-point $P$,
where $F_q$ is the $q$-th power Frobenius map and
$T_P(C)$ is the embedded tangent line at $P$ to $C$.
Needless to say, a $q$-Frobenius classical curve is
a curve which is not $q$-Frobenius nonclassical.
St\"{o}hr and Voloch \cite{sto-vol} showed that
if $C$ is $q$-Frobenius classical of degree $d$, then
\begin{equation}\label{classical}
 N_q(C) \leq \frac{1}{2}d(d+q-1),
\end{equation}
and Hefez and Voloch \cite{hef-vol} proved that
if $C$ is $q$-Frobenius nonclassical of degree $d$, then
$d \geq \sqrt{q}+1$ and
\begin{equation}\label{nonclassical}
 N_q(C) = d(q-d + 2).
\end{equation}
Each of these two estimates for $N_q(C)$ is stronger than the expected bound
if $2 \leq d \leq q-1$ for (\ref{classical})
or $d \geq \sqrt{q}+1$ for (\ref{nonclassical}).
In fact,
\[
(d-1)q + 1 - \frac{1}{2}d(d+q-1) =\frac{1}{2}(d-2)(q-d-1)
\]
and
\[
(d-1)q + 1 -d(q-d+2) = (d-\sqrt{q}-1)(d+\sqrt{q}-1).
\]
\qed

\begin{remark}\label{incorrectremark}
The St\"{o}hr-Voloch bound (\ref{classical})
is effective even if an irreducible $q$-Frobenius classical curve
$C$ has singularities.
By tracing the proof of \cite[Theorem~1]{hef-vol} carefully,
we know that if each singular point of
an irreducible $q$-Frobenius nonclassical curve $C$ is not a cusp,
then the Hefez-Voloch bound 
$ N_q(C) \leq d(q-d + 2)$
is valid for $C$
\footnote{This statement is not correct, but the conclusion of this remark is correct. See Addendum.}.
With those bounds,
taking into account the fact that
Weil's bound holds for any irreducible plane curve $C$ of degree $d$
as
$
N_q(C) \leq q+1 +(d-1)(d-2)\sqrt{q}
$
\cite[Theorem~9.57]{hir-kor-tor},
we can weaken the assumption on $C$ of Theorem~\ref{the_d<q}
as $C$ is an irreducible curve without cusp singularities.
\end{remark}

\section{Examples}
The proof of Theorem~\ref{the_d<q}, together with
Theorem~\ref{theorem_q} and the previous result~\cite{hom-kim2},
shows the following fact also.

\begin{remark}
The possible degrees $d$ of a nonsingular curve $C$ over ${\Bbb F}_q$
with $(d-1)q + 1$ rational points are
\[q+2, q+1, q, q-1, \sqrt{q} +1  
\mbox{{\rm  \ (}when $q$ is square{\rm )}, and } 2.\]
\end{remark}
For each $d$ above,
we give curves with concrete equation that attain
the bound (\ref{guess}).
\begin{itemize}
 \item Let $d = q+2$. In this case, the bound is $q^2 + q +1$
 which is the number of ${\Bbb P}^2({\Bbb F}_q)$.
 We know all irreducible or nonsingular curves of degree $q+2$
 over ${\Bbb F}_q$ that passing through all of the points of
 ${\Bbb P}^2({\Bbb F}_q)$.
 For details, see Tallini~\cite{tal} and Homma-Kim~\cite{hom-kim}.
 \item Let $d=q+1$. In the previous paper \cite{hom-kim2},
 we presented the curve
 \[
 X^{q+1} - X^2Z^{q-1} + Y^qZ - YZ^q =0
 \]
 has $q^2 + 1$ ${\Bbb F}_q$-rational points.
 \item Let $d=q$. Consider a curve $C$ defined by
 \[
 X^{q} - XZ^{q-1} + Y^{q-1}Z - Z^q =0.
 \]
 Then it is easy to see that $C$ is nonsingular and
 \[
  C({\Bbb F}_q)
    = {\Bbb P}^2({\Bbb F}_q) \setminus
        \left(
        \{ Y=0 \} \cup \{(1, \beta , 0 ) \mid \beta \in {\Bbb F}_q \}
        \right).
 \]
 Hence $N_q(C) = q^2+q+1 - 2q = (q-1)q +1.$
 \item Let $d=q-1$. As was mentioned by Sziklai~\cite{szi},
 the curve $\alpha X^{q-1} + \beta Y^{q-1} - (\alpha + \beta)Z^{q-1} =0$
 with $\alpha \beta (\alpha + \beta)\neq 0$
 has $(q-2)q + 1$ rational points.
 This curve is nonsingular and the set of rational points is
 \[
  C({\Bbb F}_q)
    = {\Bbb P}^2({\Bbb F}_q) \setminus
    \left(
     \{ X=0 \} \cup \{ Y=0 \} \cup \{ Z=0 \}
    \right).
  \]
  \item Let $q$ be a square. Then a Hermitian curve $C$
  of degree $\sqrt{q} +1$ over ${\Bbb F}_q$ attains this bound.
  Actually, $N_q(C) = (\sqrt{q})^3 +1
                     = ((\sqrt{q} +1) -1)q +1$.
  \item Let $d=3$. For a fixed field ${\Bbb F}_q$,
  there is a nonsingular curve over ${\Bbb F}_q$ with $2q+1$
  rational points if and only if
  $q= 2$ or $3$ or $4$. For details, see Schoof~\cite{sch}.
  \item Let $d=2$. It is well-known that any nonsingular quadratic
  over ${\Bbb F}_q$ has $q+1$ rational points.
\end{itemize}

\noindent
{\bf Acknowledgements.}
We would like to thank the organizers of $F_q9$ for their hospitality,
and the referee for pointing out redundancy in the original proof of 
Proposition~\ref{noirredcompnotfq}.

\vskip 8pt

\noindent
\textbf{{\Large Addendum} (Jan. 2014)}

Remark~\ref{incorrectremark} is incorrect.
In that remark, we have asserted that if an irreducible
$q$-Frobenius nonclassical plane curve $C$ of degree $d$
has no cusp singularities, then
$N_q(C) \leq d(q-d+2).$
This is not correct.
A correct assertion is that under the same assumption above,
\begin{equation}\label{correctinequality}
N_q(C) \leq (q-1)d - (2\tilde{g} -2),
\end{equation}
where $\tilde{g}$ is the genus of the normalization $\tilde{C}$
of $C$.
More precisely, let $\varphi : \tilde{C} \to C \subset {\Bbb P}^2$
be the normalization of the curve $C$ defined by $\varphi =(1, x, y)$,
where $x$ and $y$ come from coordinate functions $X/Z$ and $Y/Z$
on ${\Bbb P}^2$.
If we choose coordinates $X, Y, Z$ of ${\Bbb P}^2$ suitably,
we may assume that $x$ is a separable element of $C$.
Then we can show that
\[
{}^{\#} \{
P \in \tilde{C} \mid \varphi(P) \in C({\Bbb F}_q) \}
  = (q-1)d - (2\tilde{g} -2)
\]
by tracing the proof of \cite[Theorem 1]{hef-vol}.
Using (\ref{correctinequality}) with Weil's bound for an irreducible
plane curve of degree $d$ over ${\Bbb F}_q$,
we can prove the inequality in Theorem~4.1 for any irreducible
plane curve without cusp singularities.
Nowadays the modified Sziklai conjecture is already settled
affirmatively in our
later work \cite{hom-kim4}.
So we need not explain this approach to the partial proof
any more.


\end{document}